\documentclass[12pt]{iopart}
\usepackage[dcucite]{harvard}

\usepackage{subfigure}
\usepackage{iopams}
\usepackage{psfrag}

\usepackage[normalem]{ulem}
\usepackage[usenames,dvipsnames]{color}

\begin{document}

\note[A Fourier-based  reconstruction formula for photoacoustic computed tomography]{A simple Fourier transform-based  reconstruction formula for photoacoustic computed tomography
 with a circular or spherical measurement geometry}


\author{Kun Wang and Mark A. Anastasio\footnote[3]{To
whom correspondence should be addressed (anastasio@wustl.edu)}}

\address{
 Department
of Biomedical Engineering, 
Washington University in St. Louis, 
St. Louis, MO, 63130 USA }

\date{\today}

\begin{abstract}
Photoacoustic computed tomography (PACT), also known as optoacoustic tomography,
is an emerging imaging modality that has great potential for a wide range of 
biomedical imaging applications. 
In this Note, we derive a hybrid  reconstruction
formula that is mathematically exact and operates on a
data function that is expressed in
the  temporal frequency and  spatial domains. 
This formula explicitly reveals new insights into
 how the  spatial frequency
components of the sought-after object function
are determined by the temporal frequency components of the 
data function measured with a circular or spherical measurement geometry in
two- and three-dimensional implementations of PACT, respectively.
The structure of the reconstruction formula is surprisingly simple compared with existing 
Fourier-domain reconstruction formulae.
It also yields a straightforward numerical implementation that
is robust and two orders of magnitude more
 computationally efficient than filtered backprojection algorithms.
\end{abstract}


\maketitle

Photoacoustic computed tomography (PACT) is an emerging imaging
modality that has great potential for a wide range of biomedical imaging 
applications \cite{OraBook,WangTutorial,KrugerMP1999}. 
In PACT,  biological tissues of interest are illuminated by
use of  short laser pulses,
 which results in the generation of internal acoustic wavefields
 via the thermoacoustic effect \cite{XuReview,xu:water}. The initial amplitudes of the
 induced acoustic wavefields are proportional to the spatially variant
 absorbed optical energy density within the tissues.
 The propagated acoustic wavefields are subsequently detected
 by use of a collection of  ultrasonic transducers that are located outside the object.
 An image reconstruction algorithm is employed to estimate the absorbed optical
 energy density within the tissue
from these data.


A variety of image reconstruction algorithms have been proposed for PACT 
\cite{Xu:planar,Finch:02,Xu:2005bp,Finch:07,Kunyansky:07,Hristova:TR08,Kunyansky:12,ANU:exact,kwave:toolbox}.
While iterative image reconstruction methods hold great value due to their ability
to incorporate accurate models
of the imaging physics and the instrument response
\cite{GPaltauf:2002,HJiang:2007,Jin:2009,OATTV09:Provost,TMI:transmodel,GuoCS:2010,Huang:SPIEBoundEnhance,xu:edema,Ntziachristos:2011,Bu:2012,Kun:PMB}, 
they can lead to long 
reconstruction times, even when accelerated by use of modern computing hardware
such as graphics processing units \cite{Kun:PMB}.  This is especially problematic
in three-dimensional (3D) implementations of PACT, in which reconstruction times can
be excessively long.   Almost all experimental studies
of PACT to date have employed analytic image reconstruction algorithms.
Even if an iterative image reconstruction algorithm is to be employed,
it is often useful to employ an analytic reconstruction algorithm 
to obtain a preliminary image that can initialize the iterative algorithm
and thereby accelerate its convergence.


Most analytic reconstruction algorithms for PACT with a spherical
measurement aperture and point-like transducers have been formulated in the form of 
filtered backprojection (FBP) algorithms.  These algorithms possess a large
computational burden, requiring $O(N^5)$ floating point operations to reconstruct
a 3D image of dimension $N^3$.  Image reconstruction algorithms based on
the time-reversal principle and finite-difference schemes
require  $O(N^4)$ operations \cite{Burgholzer:07}.
Fast reconstruction algorithms for spherical measurement apertures
that require only $O(N^3 log N)$ operations
have been proposed \cite{Kunyansky:12,ANU:exact}.  
However,
numerical implementations of these formulas require
computation of special functions and multidimensional interpolation operations
in Fourier space, which require special care to avoid degradation in
reconstructed image accuracy.
It is well-known that the temporal frequency components of the pressure
data recorded on a spherical  surface are related to the Fourier
components of the sought-after object function \cite{FourierShell:07}. 
However, to date, a simple
reconstruction algorithm based on this relationship, i.e., one that does not require series expansions involving
special functions or multi-dimensional interpolations,  has yet to be developed.

In this Note, we derive a novel  reconstruction
formula for two-dimensional (2D) and 3D PACT employing circular and
spherical measurement geometries, respectively.  The mathematical
forms of the reconstruction formulae are the same in both dimensions
and are surprisingly simple compared with existing
Fourier-domain reconstruction formulae for spherical and circular measurement geometries.
The reconstruction formulae are mathematically exact 
and describe explicitly how the  spatial frequency
components of the sought-after object function are determined by the
 temporal frequency components of the measured pressure data.
 Their discrete implementations require only discrete Fourier transform,
one-dimensional interpolation, and summation operations.
A preliminary computer-simulation study is conducted to corroborate
the validity of the reconstruction formula.

We consider the canonical PACT imaging model in which the object and
surrounding medium are assumed to
possess  homogeneous and lossless acoustic properties 
 and the object is illuminated by a laser pulse with negligible temporal width.  Point-like, unfocused, ultrasonic
transducers are assumed. We also assume that the effects of the acousto-electric
impulse responses of the transducers have been deconvolved from the measured voltage
signals so that the measured data can be interpreted as pressure signals.
The 3D problem is addressed where $p(\mathbf r,t)$ denotes the photoacoustically-induced
 pressure wavefield at location  $\mathbf r \in \mathbb{R}^3$ and time $t\ge0$.
However, the analysis and reconstruction formula that follows remains valid 
for the 2D case.
The imaging physics is described by the photoacoustic wave equation 
\cite{OraBook,WangTutorial,KrugerMP1999}: 
\begin{equation}
  \nabla^2 p(\mathbf r,t) - \frac{1}{c^2}\frac{\partial^2 p(\mathbf r,t)}{\partial^2 t}=0, 
\end{equation}
subject to the initial conditions:
\begin{equation}
  p(\mathbf r,t)\Big|_{t=0}=\frac{\beta c^2}{C_p}A(\mathbf r);
  \quad \frac{\partial p(\mathbf r,t)}{\partial t} \Big|_{t=0} = 0, 
\end{equation}
where $\nabla^2$ denotes the 3D Laplacian operator and $A(\mathbf r)$ is the 
 object function to be reconstructed that is contained within the volume $V$.
Physically, $A(\mathbf r)$  represents the distribution of absorbed optical energy density.
The constant quantities $\beta$, $c$, and $C_p$ 
denote the thermal coefficient of volume expansion,
speed-of-sound,
and the specific heat capacity of the medium at constant pressure, respectively.

Let $p(\mathbf r_s, t)$ denote the pressure data recorded
at location $\mathbf r_s\in S$ on 
a spherical surface $S$ of radius $R_S$ that encloses $V$.
The continuous form of the imaging model that relates the measurement data to object function
 can be expressed as \cite{Bcox05:Greenfunc}:
\begin{equation}\label{eqn:kwaveA}
  p(\mathbf r_s, t) = \mathcal H A\equiv
      \frac{\beta c^2}{C_p(2\pi)^3}
\int_\infty\!\! d\mathbf k\,
      \hat A(\mathbf k) \cos(ckt) e^{\hat\imath \mathbf k\cdot \mathbf r_s}, 
\end{equation}
where $\mathbf k \in \mathbb{R}^3$ is the spatial frequency vector conjugate
to $\mathbf r$,  $k\equiv\vert\mathbf k\vert$, and
$ \hat A(\mathbf k)$ is the 3D Fourier transform of $A(\mathbf r)$.  We
 adopt the
 Fourier transform convention
\numparts
\begin{equation}
\hat A (\mathbf k) = \mathcal F_3 A(\mathbf r) \equiv \int_\infty d\mathbf r 
     A(\mathbf r) e^{-\hat\imath \mathbf k\cdot\mathbf r}
\end{equation}
\begin{equation}
  A (\mathbf r) = \mathcal F^{-1}_3 \hat A(\mathbf k)  \equiv
      \frac{1}{(2\pi)^3}\int_\infty d\mathbf k \hat A (\mathbf k)e^{\hat\imath \mathbf k\cdot \mathbf r}. 
\end{equation}
\endnumparts
The imaging model in Eqn.\ (\ref{eqn:kwaveA}) can be interpreted as a mapping
 $\mathcal H:O\rightarrow D$ between
infinite dimensional vector spaces that contain the object and data functions.
We will define $O$ as the vector space
of bounded and smooth functions that are compactly supported within the volume $V$.

Let the infinite set of functions   $\{\gamma_\mu(\mathbf r)\}$, indexed by
$\mu$, represent an orthonormal
basis for $O$.
The object function $A(\mathbf r)$ can be represented as
\begin{equation}
  A(\mathbf r) =
\int_\infty d\mu \,\langle A,\gamma_\mu \rangle \gamma_\mu(\mathbf r), 
\end{equation}
where the
 inner product in $O$ is defined as 
\begin{equation}
\label{eq:ip}
  \langle A,\gamma_\mu \rangle \equiv \int_V d\mathbf r A(\mathbf r)\gamma_\mu(\mathbf r)=
\frac{1}{(2\pi)^3}\int_\infty\!\! d\mathbf k\, \hat{A}(\mathbf k) \hat{\gamma}_\mu(\mathbf k),
\end{equation}
$\hat{\gamma}_\mu(\mathbf k)=\mathcal F_3 \gamma_\mu(\mathbf r)$,
and the quantity on the right-hand side of Eqn.\ (\ref{eq:ip})
follows from fact that the Fourier transform is an isometry.
A trace identity  (see Eqn.\ (1.7) in reference \cite{Finch:02} for the 3D case
and Eqn.\ (1.16) in \cite{Finch:07} for the 2D case) 
can be employed to relate the inner products in the spaces $O$ and $D$ as:
\begin{equation}\label{eqn:traceid}
 \langle A,\gamma_\mu \rangle = \frac{2C_p^2}{R_S\beta^2 c^2}
  \int_0^\infty\!\! dt \int_{S} {d \mathbf r_s}\,
   t\, p(\mathbf r_s, t) v_\mu(\mathbf r_s, t),
\end{equation}
where
\begin{equation}\label{eqn:kwaveC}
  v_\mu(\mathbf r_s, t) =\mathcal H \gamma_\mu=  \frac{\beta c^2}{C_p (2\pi)^3}
     \int_\infty\!\!\! d\mathbf k\;
      \hat \gamma(\mathbf k) \cos(ckt) 
  e^{\hat\imath \mathbf k\cdot \mathbf r_s}, 
\end{equation}
and the  right-hand side of Eqn.\ (\ref{eqn:traceid})
defines a scaled version of the inner product in $D$.

On substitution from Eqn.\ (\ref{eqn:kwaveC}) into Eqn.\ (\ref{eqn:traceid}), one obtains
\begin{equation}
\label{eq:ag1}
  \langle A,\gamma_\mu \rangle = \frac{1}{(2\pi)^3}\int_\infty d\mathbf k \,  \hat y(\mathbf k)
\hat\gamma_\mu(\mathbf k),
\end{equation}
where  
\begin{equation}
\label{eq:yk}
  \hat y(\mathbf k) \equiv \frac{2C_p}{R_S\beta}
    \int_S d\mathbf r_s e^{\hat\imath \mathbf k\cdot\mathbf r_s}
  \int_0^\infty\!\!dt\, tp(\mathbf r_s, t) \cos(ckt).
\end{equation}
Comparison of Eqns.\ (\ref{eq:ip}) and (\ref{eq:ag1}) reveals
that $\hat{A}(\mathbf k)=\hat{y}(\mathbf k)$.
By evaluating the Fourier cosine transform that is present in the right-hand side of Eqn.\ (\ref{eq:yk}),
a reconstruction formula for determining $\hat A(\mathbf k)$ can therefore be expressed as
\begin{equation}
\label{eq:key}
  \hat A(\mathbf k) = \frac{2C_p}{R_S\beta}
    \int_S d\mathbf r_s e^{\hat\imath \mathbf k\cdot\mathbf r_s}
   \mathrm{Re}\Big\{ \mathcal F_1 \big\{ tp(\mathbf r_s ,t) \big\}(\mathbf r_s, \omega)\big|_{\omega=ck}\Big\},
\end{equation}
where $\mathcal F_1$ denotes the one-dimensional (1D) Fourier transform with respect to time $t$
and 
`$\mathrm{Re}$' denotes the operation that takes the real part of quantity in the brackets. 
Subsequently, $A(\mathbf r)$ is determined as $\mathcal F_3^{-1} \hat{A}(\mathbf k)$.

Equation (\ref{eq:key}) represents a novel reconstruction for PACT and is the key result of this Note. 
Unlike previously proposed Fourier-domain reconstruction formulae
\cite{Norton:80,Kunyansky:12,ANU:exact},
Eqn.\ (\ref{eq:key}) has a  simple form and does not involve series expansions utilizing special functions. 
The reconstruction formula reveals that the measured data $p(\mathbf r_s,t)$
determine the 3D Fourier components of the $A(\mathbf r)$ via a simple process that involves
the following four steps: (1) Compute the 1D temporal Fourier transform of
the modified data function $tp(\mathbf r_s,t)$; (2) Isolate the real-valued component
of this quantity corresponding to temporal frequency $\omega=ck$;
(3) Weight this value by the plane-wave $ e^{\hat\imath \mathbf k\cdot\mathbf r_s}$;
and (4) Sum the contributions, formed in this way,  corresponding
to every  measurement location $\mathbf r_s\in S$.
This reveals the components of $\hat{A}(\mathbf k)$ residing on a sphere
of radius $k\over c$ are determined by the 1D Fourier transform of $t p(\mathbf r_s,t)$
corresponding to temporal frequency $\omega$.  In this sense, Eqn.\ (\ref{eq:key})
can be interpreted as an implementation of the Fourier Shell Identity \cite{FourierShell:07}.
Finally, the form of Eqn.\ (\ref{eq:key}) remains unchanged in the 2D case, where 
 $\mathbf r_s, \mathbf k \in \mathbb{R}^2$
and $S$ is a circle
that encloses the object.

A discrete implementation of Eqn.\ (\ref{eq:key}) possesses 
low computational complexity and desirable numerical properties.
The 1D fast Fourier transform (FFT) can be employed to approximate the action
of $\mathcal F_1$ and only a 1D interpolation is required to determine
the value of the Fourier transformed data function corresponding to
temporal frequency $\omega=ck$, where $k$ corresponds to the magnitude of vectors
$\mathbf k$ that specify a 3D Cartesian grid.   
From the values of $\hat{A}(\mathbf k)$ determined on this grid,
the 3D FFT algorithm can be employed to estimate  values
of $A(\mathbf r)$.  
If the object is represented on a $N\times N\times N$ grid and
 the number of transducer locations and time samples are both $\mathcal O(N)$,
 the computational complexity is limited
 by the 3D FFT algorithm, i.e., $\mathcal O(N^2\log N)$ in 2D and $\mathcal O(N^3\log N)$ in 3D.

\begin{figure}[c]
  \centering
  \subfigure[]{\resizebox{2.0in}{!}{\includegraphics{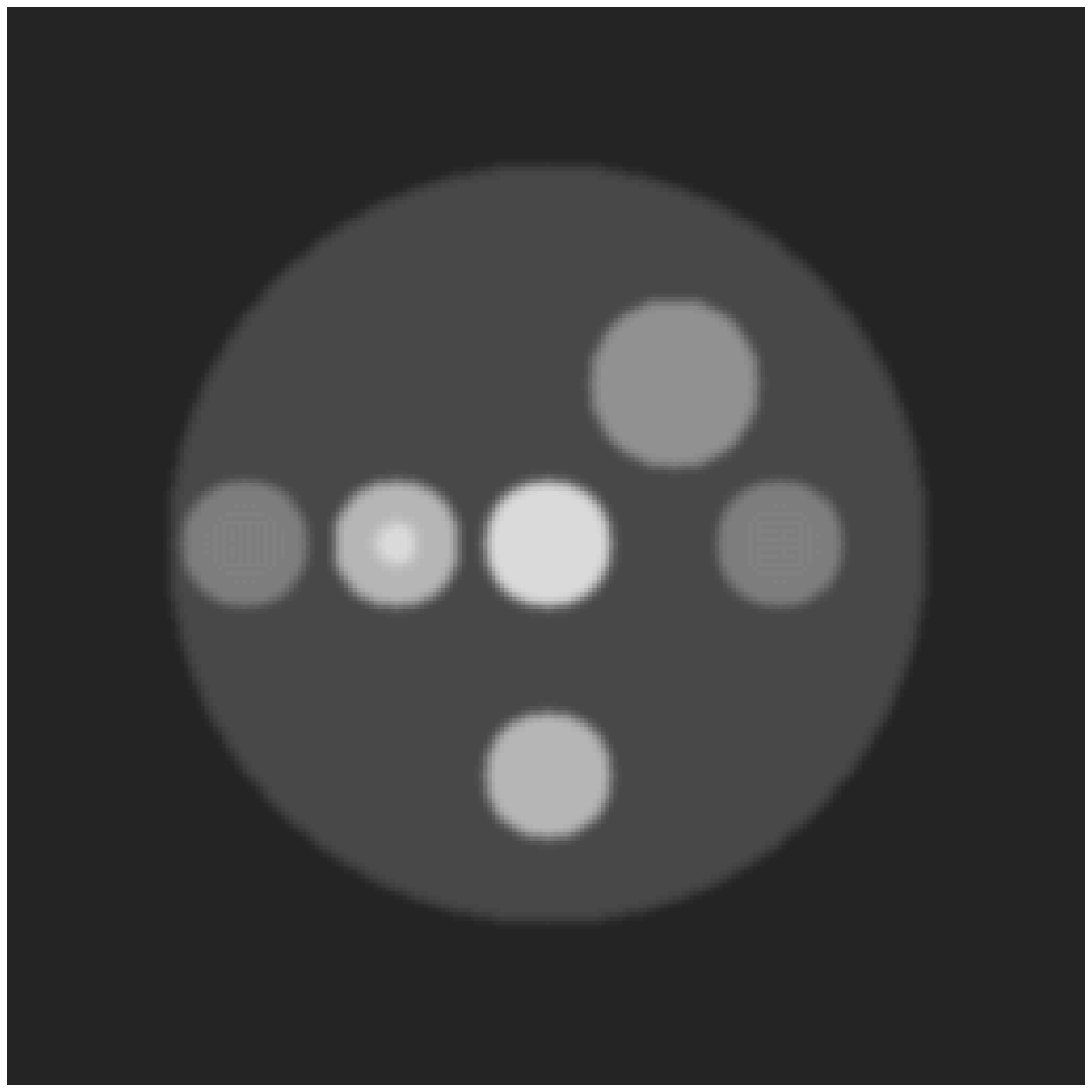}}}
  \subfigure[]{\resizebox{2.0in}{!}{\includegraphics{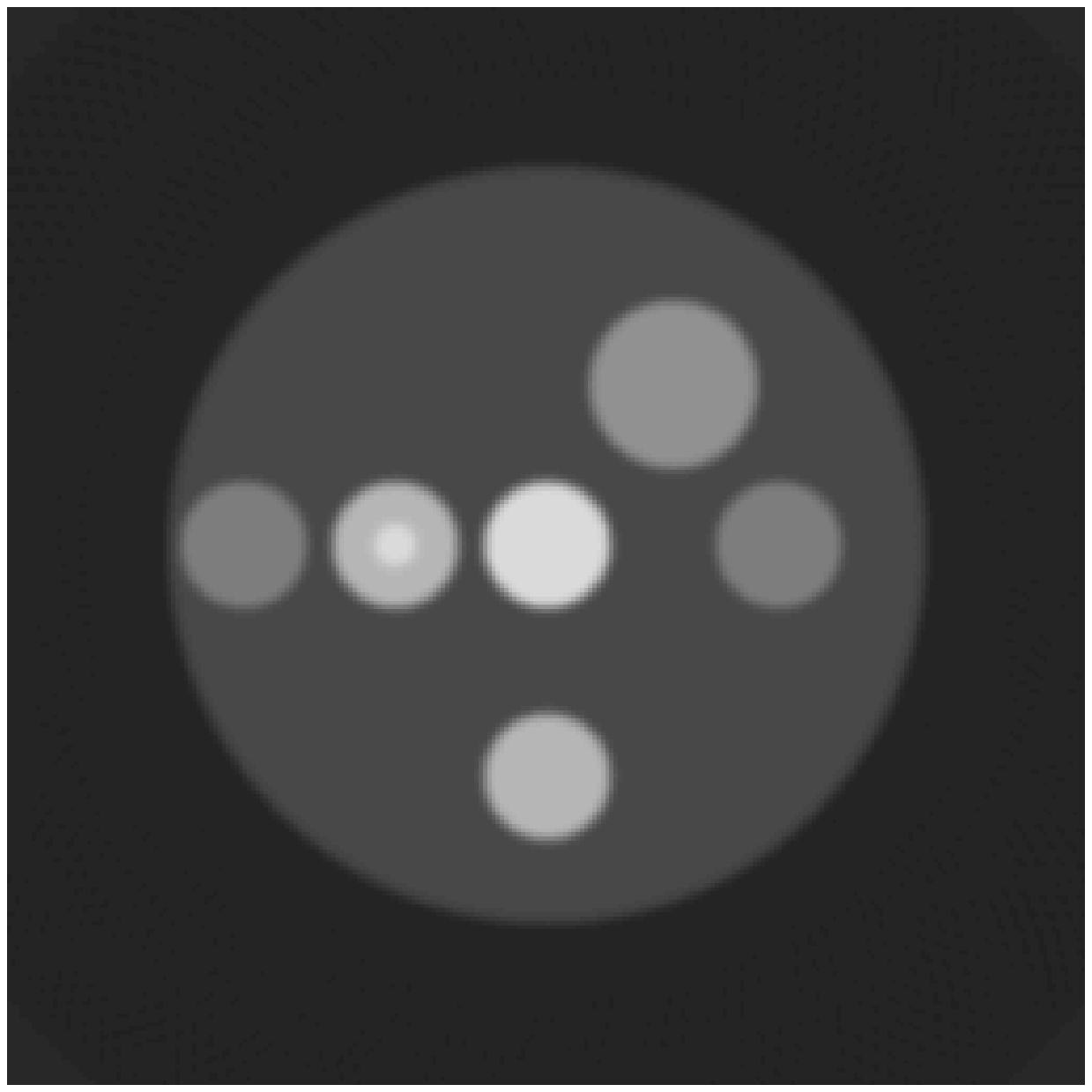}}}
  \subfigure[]{\resizebox{2.0in}{!}{\includegraphics{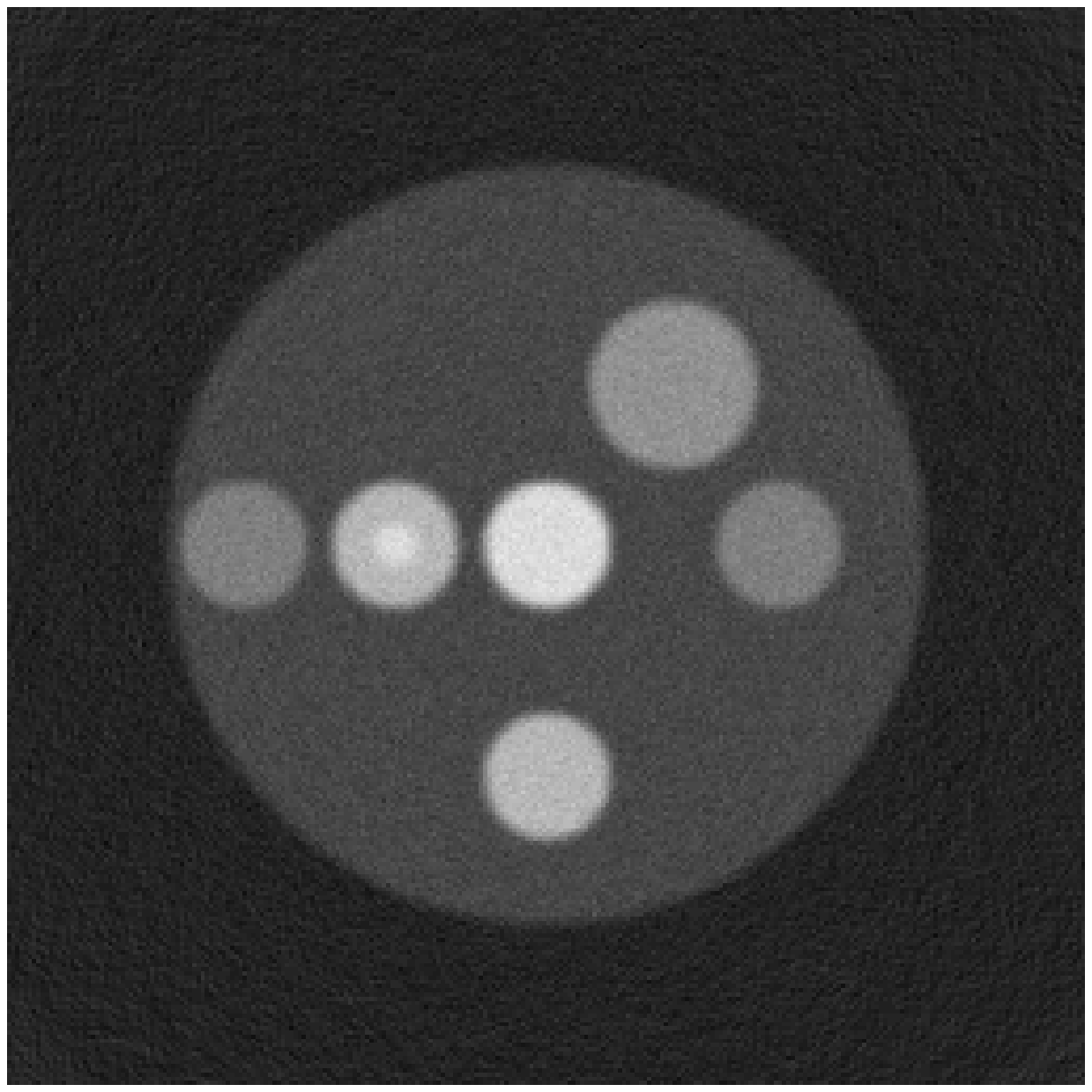}}}
  \caption{\label{fig:noiseless}
 The numerical phantom is shown in subfigure (a).
 Images reconstructed by use of 
the proposed reconstruction algorithm 
from noiseless and noisy data 
are shown in subfigures (b) and (c), respectively.
The greyscale window is $[-0.2,1.2]$. }
\end{figure}

A preliminary computer-simulation study for the 2D case was conducted to corroborate
the correctness of the reconstruction formula.
The object function $A(\mathbf r)$ was taken to be
the  numerical phantom shown in Fig.\ \ref{fig:noiseless}-(a), which
was comprised of a collection of uniform disks
that were  blurred by a 2D Gaussian kernel whose full-width-at-half-maximum was 0.3 mm.
The phantom was discretized by use of pixel expansion functions
with a pitch of 0.025 mm.
The measurement geometry consisted 
 of 256 point-like transducers that were uniformly distributed over a circle
of radius 12.8 mm that enclosed the object.
The k-wave toolbox \cite{kwave:toolbox} was employed to numerically
solve the photoacoustic wave equation and
generate simulated pressure signals at each transducer
location at a  temporal sampling rate of 30 MHz.
The simulated pressure data set generated in this way contained
$256\times 2048$ data samples.
  The speed-of-sound
and $\beta\over C_p$ were assigned values of 1.5-mm/$\mu s$
 and 1000 (arbitrary units), respectively.
A noisy data set was produced by addition of 5\% uncorrelated
Gaussian noise to the noiseless pressure data.

Images were reconstructed on a uniform 2D grid of spacing 0.1 mm 
by use of a discretized form of Eqn.\ (\ref{eq:key}) coupled with 
the 2D inverse FFT algorithm. 
In order to reconstruct images of dimension $256\times 256$, 
samples of $\hat{A}(\mathbf k)$ were determined on a uniform 2D grid of 
dimension $512\times 512$ with a sampling interval of $(0.1\times256)^{-1}$ mm$^{-1}$. 
The samples of the data function $tp(\mathbf r_s,t)$ were
zero-padded by a factor of 8 prior to estimating its 1D Fourier
transform by use of the FFT algorithm.
From these data, nearest neighbor 1D interpolation was employed to determine
the values of the term in brackets in Eqn.\ (\ref{eq:key}) corresponding to $\omega=ck$
for the sampled locations $\mathbf{k}$. 

\begin{figure}[c]
  \centering
  \subfigure[]{\resizebox{3.6in}{!}{\includegraphics{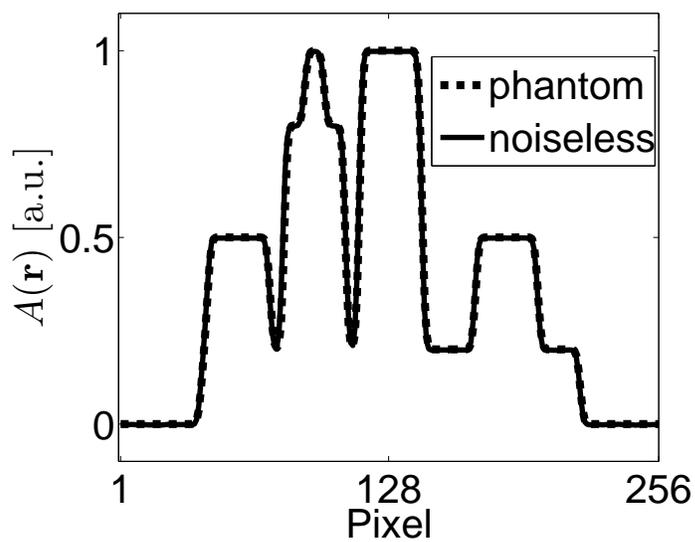}}}\\
  \subfigure[]{\resizebox{3.6in}{!}{\includegraphics{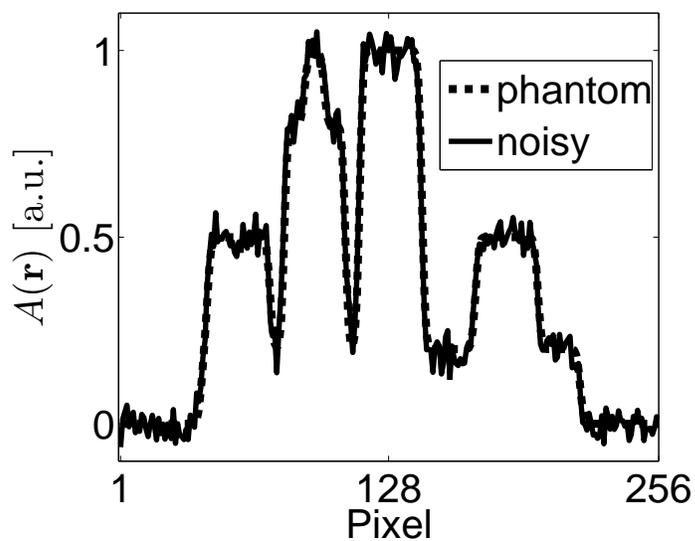}}}
  \caption{\label{fig:profile}
 Profiles corresponding to the central rows of the images
shown in Fig.\ \ref{fig:noiseless}-(b) (subfigure(a)) 
and  Fig.\ \ref{fig:noiseless}-(c) (subfigure(b)).
The solid line in subfigure (a), which corresponds to the image reconstructed
from noiseless data, almost completely
overlaps with the profile through the numerical phantom.
}
\end{figure}

The images reconstructed from the noiseless and noisy data sets
are shown in Fig.\ \ref{fig:noiseless}-(b) and (c).
 Profiles corresponding to the central rows of these images
are shown in Fig.\ \ref{fig:profile}.
These results confirm that the proposed reconstruction
algorithm can reconstruct images with high fidelity  from 
noise-free measurement data.
  Although, a systematic investigation of the noise propagation
properties of the proposed algorithm is beyond the scope of this Note,
 Figs.\ \ref{fig:noiseless}-(c) and  \ref{fig:profile}-(b)
suggest  that its performance is robust in the presence of noise.  This is to be expected, since
all operations involved in the implementation of Eqn.\ (\ref{eq:key}) are
numerically stable.

In summary, 
we have derived a Fourier-based  reconstruction
formula for PACT
employing circular and spherical measurement apertures.
 The formula is mathematically exact 
and possesses a surprisingly simple form compared with existing
Fourier-domain reconstruction formulae.
The formula yields a straightforward numerical implementation that
is  stable and is
 two orders of magnitude more
 computationally efficient than 3D filtered backprojection algorithms.
The proposed formula serves as an alternative to existing fast Fourier-based 
reconstruction formulae. 
A systematic comparison of the proposed reconstruction formula
with existing formulae by use of experimental data 
remains an important topic for future studies.


\section*{Acknowledgment}
This work was supported in part by NIH awards EB010049 and CA167446.


\section*{References}

\bibliographystyle{harvard}
\bibliography{reflect}


\clearpage

\end{document}